\documentclass{article}
\usepackage{graphicx}
\usepackage{amsmath}
\usepackage{amsfonts}
\usepackage{amssymb}

\begin{document}

\begin{center}
\textbf{Mark sequences in digraphs\bigskip}

\textbf{S. Pirzada}$^{1}$\textbf{ and U. Samee}$\bigskip^{2}$
\end{center}

$^{1}$Department of Mathematics, University of Kashmir, Srinagar, India

Email: sdpirzada@yahoo.co.in\bigskip

$^{2}$Department of Applied Mathematics, AMU, Aligarh, India \ \ \ \ \ \ \ \ \ \ \ \ \ \ \ \ \ \ \ \ \ \ \ \ \ \ \ \ \ \ \ \ \ \ \ \ \ \ \ \ \ \ \ \ \ \ \ \ \ \ \ \ \ \ \ \ \ \ \ \ \ \ \ \ \ \ \ \ \ \ \ \ \ \ \ \ \ \ \ \ \ \ \ \ \ \ \ \ \ \ \ \ \ \ \ \ \ \ \ \ \ \ \ \ \ \ \ \ \ \ \ \ \ \ \ \ \ \ \ \ \ \ \ \ \ \ \ \ \ \ \ \ \ \ \ \ \ \ \ \ \ 

\bigskip

\textbf{Abstract.} A k-digraph is an orientation of a multi-graph that is
without loops and contains at most k edges between any pair of distinct
vertices. We obtain necessary and sufficient conditions for a sequence of
non-negative integers in non-decreasing order to be a sequence of numbers,
called marks (k-scores), attached to vertices of a k-digraph. We characterize
irreducible mark sequences in k-digraphs and uniquely realizable mark
sequences in 2-digraphs.\bigskip

\textbf{1. Introduction.} Let D be a k-digraph with vertex set V = \{ v$_{1}$
,v$_{2}$ ,..., v$_{n}$\}, and let d$^{+}$(v$_{i}$) and d$^{-}$(v$_{i}$) denote
the outdegree and indegree, respectively, of a vertex v$_{i}$ . Define
p$_{v_{i}}$ (or p$_{i}$) = k(n-1)+d$^{+}$(v$_{i}$)-d$^{-}$(v$_{i}$), as the
mark of v$_{i}$, so that $0\leq p_{v_{i}}\leq2k(n-1).$ The sequence P =
[p$_{i}$]$_{1}^{n}$ in non-decreasing order is called the mark sequence of D .

A k-digraph can be interpreted as the result of a competition in which the
participants play each other at most k times, with an arc from u to v if and
only u defeats v . A player receives two points for each win, and one point
for each tie (draw), that is the case in which the two players do not play one
another or the competition between the players yields no result. With this
marking system, player v obtains a total of p$_{v}$ points.

A sequence P of non-negative integers in non- decreasing order is said to be
realizable if there exists an k-digraph with mark sequence P.

Any undefined terms are found in [3,5], and one should also take into account
of the non-standard definitions and notations introduced in this paper.

In a k-digraph, if there are x$_{1}$ arcs directed from vertex u to vertex v,
and x$_{2}$ arcs directed from vertex v to vertex u, with 0 $\leq$ x$_{1}$ ,
x$_{2}$ $\leq$ k and 0 $\leq$ x$_{1}$ + x$_{2}$ $\leq$ k, we denote it by
u(x$_{1}$-x$_{2}$)v.

We have one of the following six possibilities between any two vertices u and
v in a 2 -digraph.

(i) Exactly two arcs directed from u to v, and no arc directed from v to u,
and this is denoted by u(2 - 0)v.

(ii) Exactly two arcs directed from v to u, and no arc directed from u to v,
and this is denoted by u(0-2)v.

(iii) Exactly one arc from u to v, and exactly one arc from v to u, and this
is denoted by u(1 - 1)v.

(iv) Exactly one arc from u to v, and no arc from v to u. This is denoted by
u(1 - 0)v.

(v) Exactly one arc from v to u, and no arc from u to v, and is denoted by u(0
- 1)v.

(vi) No arcs from u to v, and no arc from v to u, and is denoted by u(0-0)v.

We note that a 1-digraph is an oriented graph, and a complete 1-digraph is a
tournament. A k-digraph D is said to be complete if there are exactly k arcs
between any pair of vertices of D.

A k-triple in a k-digraph is an induced k-subdigraph with three vertices, and
is of the form u(x$_{1}$-x$_{2}$)v(y$_{1}$-y$_{2}$)w(z$_{1}$-z$_{2}$)u, where
for i = 1 , 2 , we have 0 $\leq$ x$_{i}$ , y$_{i},$ z$_{i}$ $\leq$ k and
$0\leq\sum_{i=1}^{2}x_{i},\sum_{i=1}^{2}y_{i},\sum_{i=1}^{2}z_{i}\leq k$.
Also, in a k-digraph a 1- triple is an induced 1-subdigraph with three
vertices. A 1-triple is said to be transitive if it is of the form
u(1-0)v(1-0)w(0-1)u, or u(1-0)v(0-1)w(0-0)u, or u(1-0)v(0-0)w(0-1)u, or
u(1-0)v(0-0)w(0-0)u, or u(0-0)v(0-0)w(0-0)u, otherwise it is said to be
intransitive. A k-triple is said to be transitive if it contains only
transitive 1-triples, and a k-digraph is said to be transitive if every of its
k- triples is transitive.

A tournament is an irreflexive, complete, asymmetric digraph, and the score
s$_{v}$ of a vertex v in a tournament is the number of arcs directed away from
that vertex, and the score sequence S(T) of a tournament T is formed by
listing the vertex scores in non-decreasing order. The following criterion is
given by Landau [4].\bigskip

\textbf{Theorem 1.1 [4]. }A sequence [s$_{i}$]$_{1}^{n}$ of non-negative
integers in non-decreasing order is the score sequence of a tournament if and
only if

\begin{center}
$\sum_{i=1}^{k}s_{i}\geq\left(
\begin{array}
[c]{c}%
k\\
2
\end{array}
\right)  ,$\ \ \ \ \ \ \ for 1$\leq$ k$\leq$ n,
\end{center}

and equality for k = n.

With the marking system, the mark p$_{v}$ of a vertex v in a tournament is
given by p$_{v}$ = 2s$_{v}$ + n - 1, and Landau's conditions become

\begin{center}
$\sum_{i=1}^{k}p_{i}\geq k(n+k-2),$\ \ \ \ \ \ \ for 1$\leq$ k$\leq$ n,
\end{center}

with equality for k = n.\bigskip

An oriented graph is a digraph with no symmetric pairs of directed arcs and
without self-loops. Avery [2] defined a$_{v}$= n-1+d$^{+}$(v)-d$^{-}$(v), 0
$\leq$ a$_{v}\leq$ 2n - 2, as the score of a vertex v in an oriented graph D,
and A = [a$_{1}$, a$_{2}$,...,a$_{n}$] in non-decreasing order is the score
sequence of D. The following result is due to Avery and a constructive proof
can be found in [8].

\textbf{Theorem 1.2 [2].} A sequence A = [a$_{i}$]$_{1}^{n}$ of non-negative
integers in non-decreasing order is the score sequence of an oriented graph if
and only if

\begin{center}
$\sum_{i=1}^{k}a_{i}\geq k(k-1),$\ \ \ \ \ \ \ for 1$\leq$ k$\leq$ n,
\end{center}

with equality for k = n.

Once again, with the marking system, the mark p$_{v}$ of a vertex v in an
oriented graph is given by p$_{v}$ = a$_{v}$ + n - 1, and Avery's conditions become

\begin{center}
$\sum_{i=1}^{k}p_{i}\geq k(n+k-2),$\ \ \ \ \ \ \ for 1$\leq$ k$\leq$ n,
\end{center}

with equality for k = n.\bigskip

A k-digraph D is said to be complete if there are exactly k arcs between every
pair of vertices of D. If in a k-digraph D there are exactly k arcs, which are
parallel, between every pair of vertices of D, then D is called a k
tournament. A double tournament can be treated as a tournament whose arcs have
been duplicated.\bigskip

The following result can be easily established, and is analogues to Theorem
2.2 of Avery [2].\bigskip

\textbf{Lemma 2.1.} If D and D$^{/}$ are two k-digraphs with the same mark
sequence, then D can be transformed to D$^{/}$ by successively transforming
(i) appropriate 1-triples in one of the following ways,

either (a) by changing the intransitive 1-triple u(1-0) v(1-0) w(1-0)u to a
transitive 1-triple u(0-0) v(0-0) w(0-0)u, which has the same mark sequence,
or vice versa,

or (b) by changing an intransitive 1-triple u(1-0) v(1-0) w(0-0)u to a
transitive 1-triple u(0-0) v(0-0) w(0-1)u, which has the same mark sequence,
or vice versa.

or (c) by changing a double u(1-1)v to a double u(0-0)v which has the same
mark sequence, or vice versa.\bigskip

We note here that in a transitive tournament T, all its 1-triples are of the
form u(1- 0)v(1- 0)w(0- 1)u, for all vertices u, v and w in T. Similarly in a
transitive oriented graph, all the 1-triples are of the form u(1 - 0)v(1 -
0)w(0 - 1)u, u(1 - 0)v(0 - 1)w(0 - 0)u, u(1 - 0)v(0 - 0)w(0 - 1)u, u(1 - 0)v(0
- 0)w(0 - 0)u, u(0 - 0)v(0 - 0)w(0 - 0)u. Clearly, in the transitive double
tournament D, we have u(2-0)v(2-0)w(0-2)u for all vertices u, v and w in D.\bigskip

Now, we have the following observation.\bigskip

\textbf{Theorem 2.1. }Among all k-digraphs with a given mark sequence those
with the fewest arcs are transitive.

\textbf{Proof.} Let P be a mark sequence, and let D be a realization of P that
is not transitive. Then D contains an intransitive 1-triple. If it is of the
form u(1-0) v(1-0) w(1-0)u, it can be transformed by operation i(a) of Lemma
2.1 to a transitive 1-triple u(0-0) v(0-0) w(0-0)u with the same mark sequence
and three arcs fewer. If D contains an intransitive 1-triple of the form
u(1-0) v(1-0) w(0-0)u, it can be transformed by operation i(b) of Lemma 2.1 to
a transitive 1-triple of the form u(0-0) v(0-0) w(0-1)u with the same mark
sequence and one arc fewer. If D contains both types of intransitive
1-triples, then again they can be transformed to transitive 1-triples, and
certainly with lesser arcs. In case D contains a double u(1-1)v, it can be
transformed to u(0-0)v by operation of Lemma 2.1 with the same mark sequence
and two arcs fewer.\bigskip\ \ \ $\blacksquare$

The following result is the existence criteria for realizability of mark
sequences in k-digraphs.\bigskip

\bigskip\textbf{Theorem 2.2.} A sequence [p$_{i}$]$_{1}^{n}$ of non-negative
integers in non-decreasing order is the mark sequence of a k-digraph if and
only if

\begin{center}
$\sum_{i=1}^{t}p_{i}\geq kt(t-1),$\ \ \ \ \ \ \ for 1$\leq$ t$\leq$ n,
\end{center}

with equality for t = n.

\textbf{Proof. (i) Sufficiency.} Let q$_{i}$ = p$_{i}$-k(n-1). Then,
$\sum_{i=1}^{n}q_{i}=0$ and we may assume that $q_{1}$\bigskip$\leq q_{2}%
\leq...\leq q_{r}<0\leq q_{r+1}\leq...\leq q_{n}.$

Construct a network with vertex set \{s, v$_{1}$, v$_{2}$,..., v$_{n}$, t\} of
cardinality n+2 as follows.

1. There are arcs (s, v$_{i}$), $1\leq i\leq r$ from the source s to vertex
v$_{i}$. The arc (s, v$_{i}$) has capacity -q$_{i}$, $1\leq i\leq r.$

2. Arcs (v$_{i}$, t) from v$_{i}$ to the sink t, $r+1\leq i\leq n.$ The arc
(v$_{i}$, t) has capacity -q$_{i}.$

3. For each pair v$_{i},$ v$_{j}$ of distinct vertices (i $\neq$\ j), we have
one arc from v$_{i}$ to v$_{j}$and one arc from v$_{j}$ to v$_{i}$, each with
capacity k.

It is easy to check that a k-digraph with mark sequence [p$_{i}$]$_{i}^{n}$
can be obtained from an integral flow of value $-\sum_{i=1}^{r}q_{i}%
=\sum_{i=r+1}^{n}q_{i}$ by reducing the flow on cycles of length 2 until one
of the two edges has flow value zero.

In view of the max-flow-min-cut-Theorem, it suffice to check that each cut has
capacity at least $\sum_{i=r+1}^{n}q_{i}.$

We thus assume that \{s\}$\cup$ C is a cut, C $\subseteq$ \{v$_{1}$, v$_{2}%
$,..., v$_{n}$\}, $\left|  C\right|  $ = t, and that $\left|  C\cap
\{v_{1},v_{2},...,v_{r}\}\right|  =a,\ \left|  C\cap\{v_{r+1},v_{r+2}%
,...,v_{n}\}\right|  =b=t-a.$

For its capacity, we have the following estimate.

cap (\{s\}$\cup$ C) $=\sum_{i:i\leq r,v_{i}\notin C}-q_{i}+\sum_{i:i>r,v_{i}%
\in C}q_{i}+t(n-t)\cdot k\smallskip$

\ \ \ \ \ \ \ \ \ \ \ \ \ \ \ \ \ \ \ \ \ $\geq-\sum_{i=a+1}^{r}q_{i}%
+\sum_{i=r+1}^{r+b}q_{i}+t(n-t)\cdot k\smallskip.$

This expression is bounded from below by $-\sum_{i=1}^{r}q_{i}=\sum
_{i=r+1}^{n}q_{i}$

if and only if $\ \ \ \ \sum_{i=1}^{a}q_{i}+\sum_{i=r+1}^{r+b}q_{i}%
+t(n-t)\cdot k\geq0\smallskip$

if and only if $\ \ \ \ \sum_{i=1}^{a}p_{i}+\sum_{i=r+1}^{r+b}p_{i}%
+t(n-t)\cdot k\geq t\cdot k(n-1)$

\ \ \ \ \ \ (since p$_{i}$ = k(n-1)+q$_{i}$)

if and only if $\ \ \ \ \sum_{i=1}^{a}p_{i}+\sum_{i=r+1}^{r+b}p_{i}\geq kt(t-1).$

This latter inequality is certainly implied by the inequality

\begin{center}
$\sum_{i=1}^{t}p_{i}\geq kt(t-1),$
\end{center}

since the p$_{i}$ are non-decreasing.

\textbf{(ii) Necessity.} Follows from the construction in (i) if we use the cuts

\ \{s\}$\cup$\{v$_{1}$, v$_{2}$,..., v$_{t}$\}, \ $1\leq t\leq n.\bigskip$
\ \ \ \ \ \ \ \ \ $\blacksquare$

The following result is the existence criteria for realizability of mark
sequences in 2-digraphs.\bigskip\ The proof follows from Theorem 2.2. Here we
give a different proof.\bigskip

\bigskip\textbf{Theorem 2.3.} A sequence [p$_{i}$]$_{1}^{n}$ of non-negative
integers in non-decreasing order is the mark sequence of a 2-digraph if and
only if

$\sum_{i=1}^{k}p_{i}\geq2k(k-1),$\ \ \ \ \ \ \ for 1$\leq$ k$\leq$ n, \ \ \ \ \ \ \ \ \ \ \ \ \ (2.3.1)

with equality for k = n.

\textbf{Proof. Necessity.} Let D be a 2-digraph with mark sequence [p$_{i}%
$]$_{1}^{n}$. Let W be the 2-subdigraph induced by any set of k vertices
w$_{1}$, w$_{2}$,...,w$_{k}$ of D. Let $\alpha$ denote the number of arcs of D
that start in W and end outside W, and let \ss\ denote the number of arcs of D
that start outside of W and end in W. Note that each vertex w in W, and for
every vertex v of D not in W, there are atmost two arcs from v to w, so that
$\beta\leq$ 2k(n - k). Therefore, we have $\beta\leq$ 2nk - 2k$^{2}$. Then,
$\sum_{i=1}^{k}p_{w_{i}}=\sum_{i=1}^{k}(2n-2+d_{D}^{+}(w_{i})-d_{D}^{-}%
(w_{i}))\smallskip$

\ \ \ \ \ \ \ \ \ $=2nk-2k+\sum_{i=1}^{k}d_{D}^{+}(w_{i})-\sum_{i=1}^{k}%
d_{D}^{-}(w_{i})\smallskip$

\ \ \ \ \ \ \ \ \ \ $=2nk-2k+\left[  \sum_{i=1}^{k}d_{W}^{+}(w_{i}%
)+\alpha\right]  -\left[  \sum_{i=1}^{k}d_{W}^{-}(w_{i})\smallskip
+\beta\right]  \smallskip$

\ \ \ \ \ \ \ \ \ \ $=2nk-2k+$ (number of arcs of W) + $\alpha$ - (number of
arcs of W) - $\beta\smallskip$

\ \ \ \ \ \ \ \ \ \ $=2nk-2k+\alpha-\beta$\ \ \ \ \ \ \ \ \ \ \ \ \ \ \ \ \ \ \ \ \ \ \ \ \ \ \ \ \ \ \ \ \ \ \ \ \ \ \ \ \ \ \ \ (2.3.2)

\ \ \ \ \ \ \ \ \ \ $\geq2nk-2k-\beta\smallskip$

\ \ \ \ \ \ \ \ \ \ $\geq$\ $2nk-2k-2nk+2k^{2}=2k(k-1).$

Applying this result to the k vertices with marks p$_{1}$, p$_{2}$,...,
p$_{k}$ yields the desired inequality. If k = n, then $\alpha=\beta=0$, and
the required equality follows from Equation (2.3.2).

\textbf{Sufficiency.} This is proved by contradiction .Assume all sequences of
non-negative integers in non-decreasing order of length fewer than n,
satisfying conditions (2.3.1) be the mark sequences. Let n be the smallest and
with this choice of n, p$_{1}$ be the smallest possible such that P = [p$_{i}%
$]$_{1}^{n}$\ is not a mark sequence . Two cases arise,

(a) equality in (2.3.1) holds for some k
$<$%
n, and

(b) each inequality in (2.3.1) is strict for all k
$<$%
n.

\textbf{Case (a).} Assume k (k
$<$%
n) is the smallest such that

\begin{center}
$\sum_{i=1}^{k}p_{i}=2k(k-1).$
\end{center}

Clearly, the sequence [p$_{1}$, p$_{2}$,..., p$_{k}$] satisfies conditions
(2.3.1), and is a sequence with length less than n. So, by assumption,
[p$_{i}$]$_{1}^{k}$ is a mark sequence of some 2-digraph, say D$_{1}$.

Further,

\ \ \ $\sum_{i=1}^{m}(p_{k+i}-4k)=\sum_{i=1}^{m+k}p_{i}-\sum_{i=1}^{k}%
a_{i}-4mk\smallskip$

\ \ \ \ \ \ \ \ \ \ \ \ \ \ \ \ \ \ \ \ \ \ \ \ \ \ \ \ \ $\geq
2(m+k)(m+k-1)-2k(k-1)-4mk\smallskip$

\ \ \ \ \ \ \ \ \ \ \ \ \ \ \ \ \ \ \ \ \ \ \ \ \ \ \ \ \ $=2m(m-1),$\ \ \ \ 

for each m, 1 $\leq$ m $\leq$ n-k, with equality when m = k. As m
$<$%
n, thus by the minimality of n, the sequence [p$_{k+1}$-4k, p$_{k+2}$-4k,...,
p$_{n}$-4k] is the mark sequence of some 2-digraph D$_{2}$. The 2-digraph D of
order n consisting of disjoint copies of D$_{1}$ and D$_{2}$, such that
u(2-0)v for each vertex u $\in$ D$_{2}$ and for each vertex v $\in$ D$_{1}$,
has mark sequence P = [p$_{i}$]$_{1}^{n}$, which is a contradiction.

\textbf{Case (b).} Assume that each inequality in condition (2.3.1) is strict
for all k
$<$%
n. Obviously, p$_{1}$%
$>$%
0. Consider the sequence P$^{/}$ = [p$_{i}^{/}$]$_{1}^{n}$, defined by

\begin{center}
$p_{i}^{/}\left\{
\begin{array}
[c]{c}%
p_{i}-1,\ \ \ \ \ i=1,\\
p_{i}+1,\ \ \ \ \ i=n,\\
p_{i},\ \ \ \ \ \ otherwise.
\end{array}
\right.  $
\end{center}

Then,

\bigskip\ \ \ \ $\sum_{i=1}^{k}p_{i}^{/}=$\ $\left(  \sum_{i=1}^{k}%
p_{i}\right)  -1>2k(k-1)-1=2k(k-1),$\ 

for all k, 1 $\leq$ k
$<$%
n, and

\ \ \ \ $\sum_{i=1}^{n}p_{i}^{/}=$\ $\left(  \sum_{i=1}^{n}p_{i}\right)
-1+1=2n(n-1).$\ \ \ \ \ \ \ \ \ \ 

This shows that the sequence P$^{/}$ = [p$_{i}^{/}$]$_{1}^{n}$ satisfies
condition (2.3.1), and therefore is a mark sequence of some 2-digraph D. Let u
and v denote the vertices with mark $p_{i}^{/}=p_{i}-1$ and $p_{n}^{/}%
=p_{n}-1$ respectively.

If in D, u(0-2)v, or u(1-1)v, or u(1-0)v, or u(0-1)v, or u(0-0)v, then
transforming them respectively to u(0-1)v, or u(1-0)v, or u(2-0)v, or u(1-1)v,
or u(1-0)v, we obtain a 2-digraph with mark sequence P, a contradiction.

In D, let u(2-0)v. We have $p_{v}^{/}\geq p_{u}^{/}+2.$ If there exists at
least one vertex w $\in$ D-\{u , v\} such that the 2-triples formed by the
vertices u, v and w contain an intransitive 1-triple of the form
u(1-0)v(1-0)w(1-0)u, or u(1-0)v(1-0)w(0-0)u, or u(1-0)v(0-0)w(1-0)u,
transforming them respectively to u(1-0)v(0-0)w(0-0)u, or u(1-0)v(0-0)w(0-1)u,
or u(1-0)v(0-1)w(0-0)u, we obtain a 2-digraph with mark sequence P, which is a contradiction.

Assume for each vertex w $\in$ D-\{u , v\}, the 2-triples formed by the
vertices u, v and w contain only transitive 1-triples of the form

(i) u(1-0)v(1-0)w(0-1)u, (ii) u(1-0)v(0-1)w(1-0)u, (iii) u(1-0)v(0-1)w(0-1)u,
(iv) u(1-0)v(0-0)w(0-1)u, (v) u(1-0)v(0-1)w(0-0)u, (vi) u(1-0)v(0-0)w(0-0)u.

Then, clearly $p_{v}^{/}<p_{u}^{/}+2$, since d$_{u}^{+}>$ d$_{v}^{+}$ and
d$_{u}^{-}<$ d$_{v}^{-}$, and we get a contradiction.

If (i) appears for every vertex w $\in$ D-\{u , v\}, so that the 2-triples
formed by u, v and w is of the form u(2-0)v(2-0)w(0-1)u, then

\ \ \ \ \ \ $p_{v}^{/}=2n-2+$d$_{v}^{+}-$ d$_{v}^{-}=2n-2+2(n-2)-2=4n-8,$

and $\ p_{u}^{/}=2n-2+$d$_{u}^{+}-$ d$_{u}^{-}=2n-2+n-2+2=3n-2.$

Therefore, \ $p_{v}^{/}=p_{u}^{/}+n-6.$

For n
$<$%
8, clearly $p_{v}^{/}\leq p_{u}^{/}+1$, a contradiction.

For n $\geq$ 8, we do have $p_{v}^{/}\geq p_{u}^{/}+2$, but then
u(2-0)v(2-0)w(0-1)u can be transformed to u(2-0)v(1-0)w(0-2)u, and we get a
2-digraph with mark sequence P, a contradiction.

If (ii) appears for every vertex w $\in$ D-\{u , v\} such that the 2-triple
formed by u, v and w is of the form u(2-0)v(0-1)w(2-0)u, then

\ \ \ \ \ \ $p_{v}^{/}=2n-2+$d$_{v}^{+}-$ d$_{v}^{-}=2n-2-(n-2)-2=n-2,$

and $\ p_{u}^{/}=2n-2+$d$_{u}^{+}-$ d$_{u}^{-}=2n-2-2(n-2)=4.$

Therefore, \ $p_{v}^{/}-p_{u}^{/}=n-6,$ so that $p_{v}^{/}=p_{u}^{/}+n-6.$

For n
$<$%
8, clearly $p_{v}^{/}\leq p_{u}^{/}+1$, a contradiction.

For n $\geq$ 8, we have $p_{v}^{/}\geq p_{u}^{/}+2$. Then , transforming
u(2-0)v(0-1)w(2-0)u to u(2-0)v(0-2)w(1-0)u, we obtain a 2-digraph with mark
sequence P, again a contradiction. \ \ \ \ \ \ \ $\blacksquare\bigskip$

Some stronger inequalities on marks in 2-digraphs can be found in [7]. The
next result is the analogue of Havel-Hakimi theorem on degree sequences of
simple graph.\bigskip

\textbf{Theorem 2.4.} Let P = [p$_{i}$]$_{1}^{n}$ be a sequence of
non-negative integers in non-decreasing order, where for each i, 0 $\leq$
p$_{i}$ $\leq$ 2k(n-1). Let P$^{/}$ be obtained from P by deleting the
greatest entry p$_{n}$ (= 2k(n-1)-r, say) and (a) if r $\leq$ n-1, reducing r
greatest remaining entries by one each, or (b) if r
$>$%
n-1, reducing r-(n-1) greatest remaining entries by two each, and 2n-2-r
remaining entries by one. Then, P is a mark sequence of some k-digraph if and
only if P$^{/}$ (arranged in non-decreasing order) is a mark sequence of some k-digraph.

\textbf{Proof.} Let P$^{/}$ be a mark sequence of some k-digraph D$^{/}$. If
P$^{/}$ is obtained from P as in (a), then a k-digraph D with mark sequence P
is obtained by adding a vertex v in D$^{/}$ such that v((k-1)-0)v$_{i}$ for
those vertices v$_{i}$ in D$^{/}$ with mark v$_{i}$ = p$_{i}$-1, and
v(k-0)v$_{i}$ for those vertices v$_{i}$ in D$^{/}$ with mark v$_{i}$ =
p$_{i}$. If P$^{/}$ is obtained from P as in (b), then again a 2-digraph D
with mark sequence P is obtained by adding a vertex v in D$^{/}$ such that
v((k-1)-1)v$_{i}$ for those vertices v$_{i}$ in D$^{/}$ with mark v$_{i}$ =
p$_{i}$-2 and v((k-1)-0)v$_{i}$ for those vertices v$_{i}$ in D$^{/}$ with
mark v$_{i}$ = p$_{i}$-1.

Conversely, let P be the mark sequence of some k-digraph D. We assume D is
transitive, if not D becomes transitive by using Lemma 2.1. Let V = \{v$_{1}$,
v$_{2}$,...,v$_{n}$\} be the vertex set of D, and let p$_{n}$ = 2k(n-1)-r. If
r $\leq$ n-1, construct D such that v$_{n}$((k-1)-0)v$_{i}$ for all i, n-r
$\leq$ i $\leq$ n-1, and v$_{n}$(k-0)v$_{j}$ for all j, 1 $\leq$ j $\leq$
n-r-1. Clearly, D -v$_{n}$ realizes P$^{/}$ (arranged in non-decreasing
order). If r
$>$%
n-1, construct D such that v$_{n}$((k-1)-1)v$_{i}$ for all i, 2n-r-1 $\leq$ i
$\leq$ n-1, and v$_{n}$((k-1)-0)v$_{j}$ for all j, 1 $\leq$ j $\leq$ 2n-r-2.
Then again, D - v$_{n}$ realizes P$^{/}$ (arranged in non-decreasing order).
\ \ \ \ \ \ \ \ \ \ \ $\blacksquare\bigskip$

Theorem 2.4 provides an algorithm for determining whether a given
non-decreasing sequence P of non-negative integers is a mark sequence, and for
constructing a corresponding k-digraph. At each stage, we form P$^{/}$
according to Theorem 2.4 such that P$^{/}$ is in non-decreasing order. If
p$_{n}$ = 2k(n-1) - r, deleting p$_{n}$, and performing (a) or (b) of Theorem
2.4 according as r $\leq$ n - 1, or r
$>$%
n - 1, we get P$^{/}$ = [p$_{1}^{/}$, p$_{2}^{/}$,..., p$_{n-1}^{/}$]. If the
mark of vertex v$_{i}$ was decreased by one in this process, then the
construction yielded v$_{n}$((k-1)-0)v$_{i}$, and if it was decreased by two,
then the contradiction yielded v$_{n}$((k-1)-1)v$_{i}$. For a vertex v$_{j}$
whose mark remained unchanged, the construction yielded v$_{n}$(k-0)v$_{j}.$
If this procedure is applied recursively, then it tests whether or not P is a
mark sequence, and if P is a mark sequence, then a k-digraph with mark
sequence P is constructed.\bigskip

\textbf{Theorem 2.5.} \ Let P = [p$_{i}$]$_{1}^{n}$ be a sequence of
non-negative integers in non-decreasing order, where for each i, 0 $\leq$
p$_{i}$ $\leq$ 2k(n-1). Let P$^{/}$ be obtained from P by deleting the
greatest entry p$_{n}$ (= 2k(n-1)-r, say) and (a) if r is even, say r = 2t,
reducing t of the next greatest entries by two, or (b) if r is odd, say r =
2t+1, reducing t greatest remaining entries by two, and reducing the greatest
among the remaining entries by one. Then P is a mark sequence if and only if
P$^{/}$ (arranged in non-decreasing order) is a mark sequence.

The proof follows by using the argument as in Theorem 2.4.\bigskip

Theorem 2.5 also provides an algorithm of checking whether or not a given
non-decreasing sequence P of non-negative integers is a mark sequence and for
constructing a corresponding k-digraph. At each stage, we form P$^{/}$
according to Theorem 2.5 such that P$^{/}$ is in non-decreasing order. If
p$_{n}$ = 2k(n-1)-r, deleting p$_{n}$, and performing (a), or (b), of Theorem
2.5 according as r is even or odd, we get P$^{/}$ = [p$_{1}^{/}$, p$_{2}^{/}%
$,..., p$_{n-1}^{/}$]. If the mark of the vertex v$_{i}$ was decreased by two
in the process, then the construction yielded v$_{n}$((k-1)-1)v$_{i}$, and if
it was decreased by one, then the construction yielded v$_{n}$((k-1)-0)v$_{i}%
$. For a vertex v$_{j}$ whose mark remained unchanged, the construction
yielded v$_{n}$(k- 0)v$_{j}$. If this procedure is applied recursively, then
it tests whether or not P is a mark sequence, and if P is a mark sequence,
then a k-digraph with mark sequence P is constructed.\bigskip

\textbf{3. Irreducible mark sequences\bigskip}

A k-digraph is reducible if it is possible to partition its vertices into two
nonempty sets V$_{1}$ and V$_{2}$ in such a way that there are exactly two
arcs directed from every vertex of V$_{2}$ to each vertex of V$_{1}$, and
there is no arc from any vertex of V$_{1}$ to any vertex of V$_{2}$. If
D$_{1}$ and D$_{2}$ are k-digraphs having vertex sets V$_{1}$ and V$_{2}$
respectively, then the k-digraph D consisting of all the arcs of D$_{1}$, and
all the arcs of D$_{2}$, and exactly k arcs directed from every vertex of
D$_{2}$ to each vertex of D$_{1}$ is denoted by D = [D$_{1}$, D$_{2}$]. If
this is not possible, the k-digraph is said to be irreducible. Let D$_{1}$,
D$_{2}$,..., D$_{h}$ be irreducible k-digraphs with disjoint vertex sets. Then
D = [D$_{1}$, D$_{2}$,..., D$_{h}$] is the k-digraph having all arcs of
D$_{i}$, 1 $\leq$ i $\leq$ h, and exactly k arcs from every vertex of D$_{j}$
to each vertex of D$_{i}$, 1 $\leq$ i
$<$%
j $\leq$ h. We call D$_{1}$, D$_{2}$,..., D$_{h}$ the irreducible components
of D, and such a decomposition is called the irreducible decomposition of D. A
mark sequence P is said to be irreducible if all the k-digraphs D with mark
sequence P are irreducible.\bigskip

The following result characterizes irreducible k-digraphs.\bigskip

\textbf{Theorem 3.1.} If D is a connected k-digraph with mark sequence P =
[p$_{i}$]$_{1}^{n}$, then D is irreducible if and only if, for k = 1, 2,...,n
- 1

\ \ \ \ \ \ \ \ \ $\sum_{i=1}^{t}p_{i}>kt(t-1),$\ \ \ \ \ \ \ \ \ \ \ \ \ \ \ \ \ \ \ \ \ \ \ \ \ \ \ \ \ \ \ (3.1.1)

and $\ \sum_{i=1}^{n}p_{i}=2n(n-1).$\ \ \ \ \ \ \ \ \ \ \ \ \ \ \ \ \ \ \ \ \ \ \ \ \ \ \ \ \ \ \ (3.1.2)

\textbf{Proof. }Let D be a connected, irreducible k-digraph having mark
sequence P = [p$_{i}$]$_{1}^{n}$. Condition (3.1.2) holds, since Theorem 2.2
has already established if for any k-digraph. Condition (3.1.2) also implies
that for any integer k
$<$%
n, the k-subdigraph D$^{/}$ induced by any set of t vertices has a sum of
marks in D$^{/}$ equal to 2t(t - 1). Since D is irreducible, therefore either
there is an arc from at least one of these t vertices to at least one of the
other n-t vertices, or there is exactly one arc from at least one of the other
n-t vertices to at least one vertex in D$^{/}$. Therefore, for 1 $\leq$ t
$<$%
n-1,

$\ \ \sum_{i=1}^{k}p_{i}\geq kt(t-1)+1>kt(t-1)$

For the converse, suppose that conditions (3.1.1) and (3.1.2) hold. It follows
from Theorem 2.2 that there exists a k-digraph with mark sequence P = [p$_{i}%
$]$_{1}^{n}$. Assume such a k-digraph is reducible, and let D = [D$_{1}$,
D$_{2}$,..., D$_{h}$] be the irreducible component decomposition of D. Since
there are exactly k arcs from every vertex of D$_{j}$ to each vertex of
D$_{i}$, 1 $\leq$\ i
$<$%
\ j $\leq$\ h, D is evidently connected. If m is the number of vertices in
D$_{1}$, then m
$<$%
n, and $\ \sum_{i=1}^{m}p_{i}=km(m-1)$\ which is a contradiction to the given
hypothesis. Hence, D is irreducible. $\ \ \ \ \ \ \ \ \ \ \ \ \ \ \blacksquare
\bigskip$

We note that a disconnected k-digraph is always irreducible, since if D$_{1}$
and D$_{2}$ are the components of D, then there are no arcs between vertices
of D$_{1}$ and vertices D$_{2}$.\bigskip

The following result can be easily established.\bigskip

\textbf{Theorem 3.2. }If D is a k-digraph with mark sequence P = [p$_{i}%
$]$_{1}^{n}$, and $\sum_{i=1}^{r}p_{i}=kr(r-1),$\ $\sum_{i=1}^{t}%
p_{i}=kt(t-1),$ and \ $\sum_{i=1}^{q}p_{i}>kq(q-1),$ for r + 1 $\leq$ q $\leq$
t-1, 0 $\leq$ r
$<$%
t $\leq$ n, then the k-subdigraph induced by the vertices v$_{r+1}$, v$_{r+2}%
$,..., v$_{t}$ is an irreducible component of D with mark sequence [p$_{i}%
-kr$]$_{r+1}^{t}$.\bigskip

The mark sequence P is irreducible if D is irreducible, and the irreducible
components of P are the mark sequences of the irreducible components of D.
That is, if D = [D$_{1}$, D$_{2}$,..., D$_{h}$] is the irreducible component
decomposition of a k-digraph D with mark sequence P, then the irreducible
components P$_{i}$ of P are the mark sequences of the k-subdigraphs induced by
the vertices of D$_{i}$, 1 $\leq$ i $\leq$ h. Theorem 3.2 shows that the
irreducible components of P are determined by the successive values of k for which

\ $\sum_{i=1}^{t}p_{i}=kt(t-1),$\ \ \ \ 1 $\leq$ t $\leq$ n. \ \ \ \ \ \ \ \ \ \ \ \ \ \ \ \ \ \ \ \ \ \ \ \ \ \ \ \ \ \ \ \ (3.1.3)

This is illustrated by the following examples of 2-digraph. (i) Let P = [1, 3,
9, 12, 15, 20]. Equation (3.1.3) is satisfied for k = 2, 5, 6. Therefore,
irreducible components of P are [0], [1, 4, 7], [0] in ascending order. (ii)
Let P = [0, 5, 8, 11, 17, 19]. Here equation (3.1.3) is satisfied for k = 1,
4, 6. Therefore, irreducible components of P are [0], [1, 4, 7] and [1, 3] in
ascending order.\bigskip

A mark sequence is uniquely realizable if it belong to exactly one k-digraph.
The characterization of uniquely realizable score sequences in tournaments is
given by Avery [1], and that of oriented graphs by S.Pirzada [6]. Now as an
observation, we have the following result.\bigskip

\textbf{Theorem 3.3.} The mark sequence P of a k-digraph D is uniquely
realizable if and only if every irreducible component of P is uniquely realizable.\bigskip

The next result determines which irreducible mark sequences in 2-digraphs are
uniquely realizable.\bigskip

\textbf{Theorem 3.4.} The only irreducible mark sequences that are uniquely
realizable are [0] and [1, 3].

\textbf{Proof.} Let P be an irreducible mark sequence, and let D with vertex
set V be a 2-digraph having mark sequence P. Then D is irreducible. Therefore,
D cannot be partitioned into 2-subdigraphs D$_{1}$, D$_{2}$,...,D$_{k}$ such
that there are exactly two arcs from every vertex of D$_{\alpha}$ to each
vertex of D$_{\beta}$, 1 $\leq$ $\beta$
$<$%
$\alpha$ $\leq$ k. First assume D has n $\geq$3 vertices. Let W = \{w$_{1}$,
w$_{2}$,..., w$_{r}$\} and U = \{u$_{1}$, u$_{2}$,...,u$_{s}$\} respectively
be any two disjoint subsets of V such that r + s = n. Since D is irreducible,
(1) there do not exist exactly two arcs from every w$_{i}$(1$\leq$\ i $\leq
$\ r) to each u$_{j}$ (1 $\leq$\ j $\leq$\ s), and (2) there do not exist
exactly two arcs from every u$_{j}$ (1 $\leq$\ j $\leq$\ s) to each w$_{i}$(1
$\leq$\ i $\leq$\ s). First of all we consider case (1), and then case (2)
follows by using the same argument as in (1).

\textbf{Case (1).} There exists at least one vertex, say w$_{1}$, in W, and at
least one vertex, say u$_{1}$ in U such that either (a) w$_{1}$(1 - 1)u, or
(b) w$_{1}$(0 - 2)u$_{1}$, or (c) w$_{1}$(1 - 0)u$_{1}$, (d) w$_{1}$(0 -
1)u$_{1}$, or (e) w$_{1}$(0 - 0)u$_{1}.$

Assume w$_{i}$(2-0)u$_{j}$ for each i (1 $\leq$\ i $\leq$\ r) and j (1 $\leq
$\ j $\leq$\ s), except for i = j = 1.

If in D, either (a) w$_{1}$(1-1)u$_{1}$, or (e) w$_{1}$(0 - 0)u$_{1}$, then
transforming them respectively to w$_{1}$(0 - 0)u$_{1}$, or w$_{1}$(1 -
1)u$_{1}$, gives a 2-digraph D$^{/}$ with same mark sequence. In both cases, D
and D$^{/}$ have different number of arcs, and thus are non-isomorphic.

(b) Let w$_{1}$(0 - 2)u$_{1}$. Since there are only six possibilities between
w$_{1}$ and w$_{i}$ , therefore, for any other vertex w$_{i}$ in W we have one
of the following cases.

(i) w$_{1}$(2 - 0)w$_{i}$(2 - 0)u$_{1}$(2 - 0)w$_{1}$, (ii) w$_{1}$(1 -
1)w$_{i}$(2 - 0)u$_{1}$(2 - 0)w$_{1}$, (iii) w$_{1}$(1 - 0)w$_{i}$(2 -
0)u$_{1}$(2 - 0)w$_{1}$, (iv) w$_{1}$(0 - 1)w$_{i}$(2 - 0)u$_{1}$(2 -
0)w$_{1}$, (v) w$_{1}$(0 - 0)w$_{i}$(2 - 0)u$_{1}$(2 - 0)w$_{1}$, (vi) w$_{1}%
$(0- 2)w$_{i}$(2 - 0)u$_{1}$(2 - 0)w$_{1}$.

Transforming (i) - (v) respectively to w$_{1}$(1 - 0)w$_{i}$(1 - 0)u$_{1}$(1 -
0)w$_{1}$,w$_{1}$(0 - 1)w$_{i}$(1 - 0)u$_{1}$(1 - 0)w$_{1}$, w$_{1}$(0 -
0)w$_{i}$(1 - 0)u$_{1}$(1 - 0)w$_{1}$, w$_{1}$(0 - 2)w$_{i}$(1 - 0)u$_{1}$(1 -
0)w$_{1}$, w$_{1}$(0 - 1)w$_{i}$(1 - 0)u$_{1}$(1 - 0)w$_{1}$, gives a
2-digraph with the same mark sequence. In all these five cases, D and D$^{/}$
have different number of arcs, and thus are non-isomorphic.

If (vi) occurs in D, and also w$_{q}$(2 - 0)w$_{i}$ for 1 $\leq$ i
$<$%
q $\leq$ r, then the 2-digraph D is reducible with irreducible components
D$_{1}$, D$_{2}$,...,D$_{r}$ respectively having vertex sets V$_{1}$ =
\{u$_{1}$, u$_{2}$,...,u$_{s}$, w$_{1}$\}, V$_{2}$ = \{w$_{2}$\}, V$_{3}$ =
\{w$_{3}$\},...,V$_{k}$ =\{w$_{r}$\}.

Also for any vertex u$_{j}$ in U, since there are only six possibilities
between u$_{1}$ and u$_{j}$ , we have one of the following cases.

(vii) w$_{1}$(0 - 2)u$_{1}$(0 - 2)u$_{j}$(0- 2)w$_{1}$, (viii) w$_{1}$(0 -
2)u$_{1}$(1 - 1)u$_{j}$(0 - 2)w$_{1}$, (ix) w$_{1}$(0 - 2)u$_{1}$(1 -
0)u$_{j}$(0 - 2)w$_{1}$, (x) w$_{1}$(0 - 2)u$_{1}$(0 - 1)u$_{j}$(0 - 2)w$_{1}%
$, (xi) w$_{1}$(0 - 2)u$_{1}$(0 - 0)u$_{j}$(0 - 2)w$_{1}$, (xii) w$_{1}$(0 -
2)u$_{1}$(2- 0)u$_{j}$(0 - 2)w$_{1}$.

If any one of (vii) -(xi) appears in D, then making respectively the
transformations w$_{1}$(0-1)u$_{1}$(0-1)u$_{j}$(0-1)w$_{1}$, w$_{1}%
$(0-1)u$_{1}$(1-0)u$_{j}$(0-1)w$_{1}$, w$_{1}$(0-1)u$_{1}$(2 - 0)u$_{j}$(0 -
1)w$_{1}$, w$_{1}$(0 - 1)u$_{1}$(1 - 1)u$_{j}$(0 - 1)w$_{1}$, w$_{1}$(0 -
1)u$_{1}$(1 - 0)u$_{j}$(0 - 1)w$_{1}$, we get a 2-digraph with the same mark
sequence, but number of arcs in D and D$^{/}$ are different, and thus D and
D$^{/}$ are non-isomorphic.

If (xii) and any of (i) - (v) appear simultaneously, then there exists a
2-digraph D$^{/}$ with the same mark sequence, but D and D$^{/}$ have
different number of arcs. Thus, D and D$^{/}$ are non-isomorphic.

If (vi) and (xii) appear simultaneously, and also w$_{q}$(2 - 0)w$_{i}$ for
all 1 $\leq$ i
$<$%
q $\leq$ r, then D is reducible with the irreducible components D$_{1}$,
D$_{2}$,..., D$_{r}$ having vertex sets V$_{1}$ = \{u$_{1}$, u$_{2}$,...,
u$_{s}$, w$_{1}$\}, V$_{2}$ = \{w$_{2}$\}, V$_{3}$ = \{w$_{3}$\},...,V$_{r}$ =
\{w$_{r}$\} respectively.

(c) Let w$_{1}$(1 - 0)u$_{1}$. For any vertex w$_{i}$ in W, since there are
only six possibilities between w$_{1}$ and w$_{i}$ , we have one of the
following cases.

(i) w$_{1}$(2 - 0)w$_{i}$(2 - 0)u$_{1}$(0 - 1)w$_{1}$, (ii) w$_{1}$(1 -
1)w$_{i}$(2 - 0)u$_{1}$(0 - 1)w$_{1}$, (iii) w$_{1}$(1 - 0)w$_{i}$(2 -
0)u$_{1}$(0 - 1)w$_{1}$, (iv) w$_{1}$(0 - 1)w$_{i}$(2 - 0)u$_{1}$(0 -
1)w$_{1}$, (v) w$_{1}$(0 - 0)w$_{i}$(2 - 0)u$_{1}$(0- 1)w$_{1}$, (vi) w$_{1}%
$(0 - 2)w$_{i}$ (2 - 0)u$_{1}$(0 - 1)w$_{1}$.

For (i) - (v) making respectively the transformations

w$_{1}$(1 - 0)w$_{i}$(1 - 0)u$_{1}$(0 - 2)w$_{1}$, w$_{1}$(0 - 1)w$_{i}$(1 -
0)u$_{1}$(0 - 2)w$_{1}$, w$_{1}$(0 - 1)w$_{i}$(1 - 0)u$_{1}$(0 - 2)w$_{1}$,
w$_{1}$(1 - 1)w$_{i}$(1 - 0)u$_{1}$(2 - 0)w$_{1}$, w$_{1}$(0 - 1)w$_{i}$(1 -
0)u$_{1}$(2 - 1)w$_{1}$, we obtain a 2-digraph D$^{/}$ with the same mark
sequence, but the number of arcs in D and D$^{/}$ is not equal. Thus, D and
D$^{/}$ are non-isomorphic.

Now, for any other vertex u$_{j}$ in U, there are only six possibilities
between u$_{1}$ and u$_{j}$ , and we have one of the following cases.

(vii) w$_{1}$(1 - 0)u$_{1}$(0 - 2)u$_{j}$(0 - 2)w$_{1}$, (viii) w$_{1}$(1 -
0)u$_{1}$(1 - 1)u$_{j}$(0 - 2)w$_{1}$, (ix) w$_{1}$(1 - 0)u$_{1}$(1 -
0)u$_{j}$(0 - 2)w$_{1}$, (x) w$_{1}$(1 - 0)u$_{1}$(0 - 1)u$_{j}$(0 - 2)w$_{1}%
$, (xi) w$_{1}$(1 - 0)u$_{1}$(0 - 0)u$_{j}$(0 - 2)w$_{1}$, (xii) w$_{1}$(1 -
0)u$_{1}$(2 - 0)u$_{j}$(0 - 2)w$_{1}$.

If any one of (vii) - (xi) appears, then making respectively the transformations

w$_{1}$(2-0)u$_{1}$(0-1)u$_{j}$(0-1)w$_{1}$, w$_{1}$(2-0)u$_{1}$(1-0)u$_{j}%
$(0-1)w$_{1}$, w$_{1}$(2- 0)u$_{1}$ (2 - 0)u$_{j}$(0 - 1)w$_{1}$, w$_{1}$(2 -
0)u$_{1}$(1 - 1)u$_{j}$(0 - 1)w$_{1}$, w$_{1}$(2 - 0)u$_{1}$(1 - 0)u$_{j}$(0 -
1)w$_{1}$, we get a 2-digraph D$^{/}$ with the same mark sequence, but D and
D$^{/}$ have different number of arcs. Thus, D and D$^{/}$ are non-isomorphic.

If (xii) and one of (i)-(v) appears simultaneously, we once again arrive to
the conclusion that there exists a 2-digraph D$^{/}$ with the mark sequence P,
but D and D$^{/}$ are non-isomorphic.

Thus, we are left with the case when (vi) and (xii) appear simultaneously, and
also w$_{q}$(2 - 0)w$_{i}$ for all 1 $\leq$ i
$<$%
q $\leq$ r. But, then D is reducible having the irreducible components D$_{1}%
$, D$_{2}$,...,D$_{r}$ with vertex sets V$_{1}$ = \{u$_{1}$, u$_{2}%
$,...,u$_{s}$, w$_{1}$\}, V$_{2}$ = \{w$_{2}$\},...,V$_{r}$ = \{w$_{r}$\} respectively.

(d) Let w$_{1}$(0- 1)u$_{1}$. Since there are only six possibilities between
w$_{1}$ and w$_{i}$ , therefore for any other vertex w$_{i}$ in W, we have one
of the following cases.

(i) w$_{1}$(2 - 1)w$_{i}$(2 - 0)u$_{1}$(1 - 0)w$_{1}$, (ii) w$_{1}$(1 -
1)w$_{i}$(2 - 0)u$_{1}$(1 - 0)w$_{1}$, (iii) w$_{1}$(1 - 0)w$_{i}$(2 -
0)u$_{1}$(1 - 0)w$_{1}$, (iv) w$_{1}$(0 - 1)w$_{i}$(2 - 0)u$_{1}$(1 -
0)w$_{1}$, (v) w$_{1}$(0 - 0)w$_{i}$(2 - 0)u$_{1}$(1 - 0)w$_{1}$, (vi) w$_{1}%
$(0 - 2)w$_{i}$(2 - 0)u$_{1}$(1 - 0)w$_{1}$.

If any one of (i) - (v) appears, then making respectively the transformation
w$_{1}$(1-0)w$_{i}$(1-0)u$_{1}$(0- 0)w$_{1}$, w$_{1}$(0-1)w$_{i}$(1-0)u$_{1}%
$(0-0)w$_{1}$,w$_{1}$(0-0)w$_{i}$(1 - 0)u$_{1}$(0 - 0)w$_{1}$, w$_{1}$(0 -
2)w$_{i}$(1 - 0)u$_{1}$(0 - 0)w$_{1}$, w$_{1}$(0 - 1)w$_{i}$(1 - 0)u$_{1}$(0 -
0)w$_{1}$, gives a 2-digraph D$^{/}$ with the same mark sequence, but number
of arcs in D and D$^{/}$ is different so that D and D$^{/}$ are non-isomorphic.

If (vi) appears in D, and also if w$_{q}$(2 - 0)w$_{i}$ for all 1 $\leq$ i
$<$%
q $\leq$ r, then D becomes reducible.

Now, for any other vertex u$_{j}$ in U, there are only six possibilities
between u$_{1}$ and u$_{j}$ , and we have one of the following cases.

(vii) w$_{1}$(0 - 1)u$_{1}$(0 - 2)u$_{j}$(0 - 2)w$_{1}$, (viii) w$_{1}$(0 -
1)u$_{1}$(1 - 1)u$_{j}$(0 - 2)w$_{1}$, (ix) w$_{1}$(0 - 1)u$_{1}$(1 -
0)u$_{j}$(0 - 2)w$_{1}$, (x) w$_{1}$(0 - 1)u$_{1}$(0 - 1)u$_{j}$(0 - 2)w$_{1}%
$, (ix) w$_{1}$(0 - 1)u$_{1}$(0 - 0)u$_{j}$(0 - 2)w$_{1}$, (xii) w$_{1}$(0 -
1)u$_{1}$(2 - 0)u$_{j}$(0 - 2)w$_{1}$.

If any one of (vii) -(xi) appears in D, then making respectively the transformations

w$_{1}$(0 - 0)u$_{1}$(0 - 1)u$_{j}$(0 - 1)w$_{1}$, w$_{1}$(0 - 0)u$_{1}$(1 -
0)u$_{j}$(0 - 1)w$_{1}$, w$_{1}$(0 - 0)u$_{1}$(2 - 0)u$_{j}$(0 - 1)w$_{1}$,
w$_{1}$(0 - 0)u$_{1}$(0 - 0)u$_{j}$(0 - 1)w$_{1}$, w$_{1}$(0-0)u$_{1}%
$(1-0)u$_{j}$(0-1)w$_{1}$, gives a 2-digraph D$^{/}$ with the same mark
sequence, but number of arcs in D and D$^{/}$ is different so that D is not
isomorphic to D$^{/}$.

If (xii) and any one of (i) - (v) appear simultaneously, then once again there
exists a 2-digraph D$^{/}$ with the same mark sequence, but D and D$^{/}$ have
different number of arcs so that D and D$^{/}$ are non-isomorphic.

If (vi) and (xii) appear simultaneously, and also w$_{q}$(2 - 0)w$_{i}$ for
all 1 $\leq$ i
$<$%
q $\leq$ r, then D is reducible.

Now, let D have exactly two vertices say u and v. The only irreducible mark
sequences realizing D are [2, 2], and [1, 3]. Obviously the sequence [2, 2]
has two non-isomorphic realizations namely u(0 - 0)v and u(1 - 1)v, and [1, 3]
has the unique realization u(0 - 1)v. Thus P = [1, 3] is uniquely realizable.

If D has only one vertex, then P = [0], which evidently is uniquely
realizable. \ \ \ \ \ \ \ \ \ \ \ $\blacksquare\bigskip$

Combining Theorem 3.3 and Theorem 3.4, we have the following result in 2-digraphs.\bigskip

\textbf{Theorem 3.5.} The mark sequence P of a 2-digraph is uniquely
realizable if and only if every irreducible component of P is of the form [0]
or [1, 3].\bigskip

We observe that in the mark sequence P = [4i-4]$_{1}^{n}$ every irreducible
component is [0], and thus P is uniquely realizable. We note that mark
sequences of tournaments are not uniquely realizable. To see this, consider
the mark sequence P = [2, 4, 6] realizing the tournament T. The other
2-digraph D realized by P has vertex set \{v$_{1}$, v$_{2}$, v$_{3}$\} with
v$_{1}$(0-0)v$_{2}$(0-0)v$_{3}$(2-0)v$_{1}$.

However, we observe that a mark sequence of a tournament T is uniquely
realizable if and only if the mark sequence of the double tournament of T is
uniquely realizable.\bigskip

Now, we have the following generalization of Theorem 3.5, and the proof
follows by using the argument as in Theorem 3.5.\bigskip

\textbf{Theorem 3.6.} The mark sequence P of a k-digraph is uniquely
realizable if and only if every irreducible component of P is of the form [0]
or [1, 2k-1].\bigskip

\textbf{Acknowledgements.} The authors are thankful to the ananomious referee
for his valuable suggestions and for providing an elegant proof of Theorem
2.2, which improved the presentation of the paper.\bigskip

\bigskip\ 

\textbf{References\bigskip}

[1] P.Avery, Condition for a tournament score sequence to be simple, J. Graph
Theory, No. 4 (1980) 157-164.

[2] P.Avery, Score sequences of oriented graphs, J.Graph Theory, Vol. 15, No.
3 (1991) 251-157.

[3] F.Harary, R.Z.Norman and D.Cartwright, Structural Models: An Introduction
to the Theory of Directed Graphs, John Wiley and Sons, New York (1965).

[4] M.G.Landau, On dominance relations and the structure of animal societies:
I I I. The condition for a score structure, Bull. Math. Biophys. 15 (1953) 143-148.

[5] J.W.Moon, Topics on Tournaments, Holt, Rinehart and Winston, New York (1968).

[6] S.Pirzada, Simple score sequences in oriented graphs, Novi Sad J.
Mathematics, Vol. 33, No. 1 (2003) 25-29.

[7] S.Pirzada and T.A.Naikoo, Inequalities for marks in digraphs, Mathematical
Inequalities and Applications, Vol. 9, No. 2 (2006) 189-198.

[8] S.Pirzada, New proof of a Theorem on oriented graph scores, To appear.
\end{document}